\documentclass[a4paper]{amsart}

\usepackage{amssymb,latexsym}
\usepackage{amsmath}
\usepackage{amsfonts}
\usepackage{amsthm}
\usepackage{amssymb}

\theoremstyle{plain}

\theoremstyle{definition}

\title[High rank curves induced by rational Diophantine triples]{High rank elliptic curves induced by rational Diophantine triples}
\begin{document}


\author[A. Dujella]{Andrej Dujella}
\address{
Department of Mathematics\\
Faculty of Science\\
University of Zagreb\\
Bijeni{\v c}ka cesta 30, 10000 Zagreb, Croatia
}
\email[A. Dujella]{duje@math.hr}

\author[J. C. Peral]{Juan Carlos Peral}
\address{
Departamento de Matem\'aticas\\
Universidad del Pa\'{\i}s Vasco\\
Aptdo. 644, 48080 Bilbao, Spain
}
\email[J. C. Peral]{juancarlos.peral@ehu.es}


\maketitle

\begin{center}
{\it

{{\it\small Dedicated to the memory of our friend and coauthor Juli\'an Aguirre}}

\vspace{\baselineskip} }
\end{center}

\begin{abstract}
A rational Diophantine triple is a set of three nonzero rational $a,b,c$
with the property that $ab+1$, $ac+1$, $bc+1$ are perfect squares.
We say that the elliptic curve $y^2 = (ax+1)(bx+1)(cx+1)$ is induced by the triple $\{a,b,c\}$.
In this paper, we describe a new method for construction of elliptic curves over $\mathbb{Q}$
with reasonably high rank based on a parametrization of rational Diophantine triples.
In particular, we construct an elliptic curve induced by a rational Diophantine triple
with rank equal to $12$, and an infinite family of such curves with rank $\geq 7$,
which are both the current records for that kind of curves.
\end{abstract}

\section{Introduction}

A set $\{a_1, a_2, \dots, a_m\}$ of $m$ distinct nonzero rationals
is called {\it a rational Diophantine $m$-tuple} if $a_i a_j+1$ is a perfect square
for all $1\leq i < j\leq m$.
The first rational Diophantine quadruple
$\{\frac{1}{16}, \frac{33}{16}, \frac{17}{4}, \frac{105}{16} \}$ was found by Diophantus,
while the first Diophantine quadruple in integers $\{1, 3, 8, 120\}$
was found by Fermat. In 1969, Baker and Davenport \cite{B-D} proved that Fermat's set
cannot be extended to a Diophantine quintuple in integers.
It was proved in \cite{D-crelle} that there does not exist a Diophantine sextuple in integers
and there are only finitely many Diophantine quintuples in integers.
Recently, He, Togb\'e and Ziegler proved that there are no Diophantine quintuples
in integers \cite{HTZ}.
Euler proved that there are infinitely many rational Diophantine quintuples.
In particular, he extended Fermat's quadruple by the fifth positive rational number
$\frac{777480}{8288641}$. In 2019, Stoll \cite{Stoll}
proved that extension of Fermat's set to a rational quintuple with the same property is unique.
The first example of a rational Diophantine sextuple, the
set $\{ \frac{11}{192}, \frac{35}{192}, \frac{155}{27}, \frac{512}{27}, \frac{1235}{48}, \frac{180873}{16}\}$,
was found by Gibbs \cite{Gibbs1},
while Dujella, Kazalicki, Miki\'c and Szikszai \cite{DKMS} recently proved that there are infinitely
many rational Diophantine sextuples (see also \cite{Duje-Matija,DKP,DKP-split}).
For an overview of results on
Diophantine $m$-tuples and its generalizations see \cite{Duje-Notices}.

Let $\{a,b,c\}$ be a rational Diophantine triple. Then there exist nonnegative rationals
$r,s,t$ such that $ab+1=r^2$, $ac+1=s^2$ and $bc+1=t^2$.
In order to extend the triple $\{a,b,c\}$ to a quadruple,
we have to solve the system of equations
\begin{equation} \label{eq2}
ax+1= \square, \quad bx+1=\square, \quad cx+1=\square.
\end{equation}
We assign the following elliptic curve to the system \eqref{eq2}:
\begin{equation} \label{e1}
E:\qquad  y^2=(ax+1)(bx+1)(cx+1).
\end{equation}
We say that the elliptic curve $E$ is induced by the rational Diophantine triple $\{a,b,c\}$.


Elliptic curves induced by rational Diophantine triples were used for the first time
in the construction of elliptic curves with relatively large rank in \cite{D-rocky}
(let us mention that in \cite{HPZ} all $S$-integral points on some elliptic curves associated
with the quintuple $\{ \frac{1}{16}, \frac{33}{16}, \frac{105}{16}, 20, 1140 \}$ were computed,
which was a motivation for considering connections between elliptic curves and Diophantine $m$-tuples).
By using subtriples of certain rational Diophantine quintuples,
elliptic curves with rank $7$ over $\mathbb{Q}$ and rank $4$ over $\mathbb{Q}(t)$ were constructed
in \cite{D-rocky}.
That result was improved in \cite{D-glasnik} where several examples of curves with rank $9$ were found
by considering subtriples of the following generalization of Fermat's quadruple:
$\{k-1,k+1,4k,16k^3-4k\}$.
These results were further improved in our joint paper with Juli\'an Aguirre \cite{ADP},
where we constructed an elliptic curve with rank $11$ over $\mathbb{Q}$
(induced that the triple
$\{ \frac{795025}{3128544}, -\frac{22247424}{7791245}, \frac{24807390285149}{97501011189120} \}$)
and rank $5$ over $\mathbb{Q}(t)$. The construction was based on subtriples of quadruples
of the form $\{a, a(k + 1)^2 - 2 k, a(2 k + 1)^2 - 8 k - 4, a k^2 - 2 k - 2 \}$.
We used similar method in \cite{D-Peral-RACSAM} and constructed several new elliptic curves
with rank $11$ over $\mathbb{Q}$ and rank $6$ over $\mathbb{Q}(t)$ (see also \cite{D-Peral-JGEA}).

Note that in all mentioned results the elliptic curves have torsion group
$\mathbb{Z}/2\mathbb{Z} \times \mathbb{Z}/2\mathbb{Z}$.
The application of elliptic curves induced by rational Diophantine triples in construction of high rank
curves appears to be even more fruitful in the case of larger torsion groups.
Such curves were used in \cite{D-Peral-LMSJCM,D-Peral-JGEA} for finding elliptic curves with
the largest known rank over $\mathbb{Q}$ (rank $9$;
induced by the triples $\{ \frac{301273}{556614}, -\frac{556614}{301273}, -\frac{535707232}{290125899} \}$
and $\{ \frac{181800}{127673}, -\frac{127673}{181800}, -\frac{996869751703}{2072406375000} \}$)
and $\mathbb{Q}(t)$ (rank $4$)
with torsion group $\mathbb{Z}/2\mathbb{Z} \times \mathbb{Z}/4\mathbb{Z}$.
This construction uses triples of the form $\{a, -\frac{1}{a}, c\}$ which induce elliptic curves
with points of order $4$. It is shown in \cite{D-Peral-RACSAM} that the elliptic curve
with largest known rank over $\mathbb{Q}$ (rank $6$; originally found by Elkies in 2006)
with torsion group $\mathbb{Z}/2\mathbb{Z} \times \mathbb{Z}/6\mathbb{Z}$ is induced by the triple
$\{ \frac{31269599}{31628160}, -\frac{23721120}{31269599}, \frac{1461969791}{7144352640} \}$.

Furthermore, it was shown in \cite{D-glasnik} that every elliptic curve with torsion
group $\mathbb{Z}/2\mathbb{Z} \times \mathbb{Z}/8\mathbb{Z}$ is induced by a Diophantine triple
(see also \cite{CG}). In particular, the triple
$\{ \frac{408}{145}, -\frac{145}{408}, -\frac{145439}{59160} \}$ induces the curve with
torsion group $\mathbb{Z}/2\mathbb{Z} \times \mathbb{Z}/8\mathbb{Z}$ and rank $3$ over $\mathbb{Q}$,
found by Connell and Dujella in 2000, what is the largest known rank for curves with that torsion group.

Although in the case of torsion group $\mathbb{Z}/2\mathbb{Z} \times \mathbb{Z}/2\mathbb{Z}$,
the record ranks over $\mathbb{Q}$ (rank $15$) and $\mathbb{Q}(t)$ (rank $7$)
were discovered by Elkies \cite{Elkies1,Elkies2} with different methods,
we believe that it is still interesting question to investigate how large can be the rank
of elliptic curves induced by rational Diophantine triples.
In this paper, we construct an elliptic curve induced by a rational Diophantine triples
with rank equal to $12$, and an infinite family of such curves with rank $\geq 7$,
which both improve previous results of the type.

\section{Construction of an  elliptic curve with rank $12$}

By the coordinate transformation $x\mapsto \frac{x}{abc}$, $y\mapsto \frac{y}{abc}$, applied
to the curve $E$, we obtain the equivalent curve
\begin{equation} \label{e2}
E':\qquad  y^2=(x+ab)(x+ac)(x+bc).
\end{equation}
The curve $E'$ has three $2$-rational points $A=[-bc,0]$, $B=[-ac,0]$, $C=[-ab,0]$,
and other two rational points $P=[0,abc]$ and $S=[1,rst]$, where
$ab+1=r^2$, $ac+1=s^2$, $bc+1=t^2$. We may expect that in general the points $P$ and $S$
will be independent points of infinite order, so that the rank of $E'$ will be at least $2$.

To increase the rank, we will use the parametrization of rational Diophantine triples
due to Lasi\'c \cite{Lasic}
(see also \cite{DKP-split}):
\begin{align*}
a &= \frac{2 t_1 (1 + t_1 t_2 (1 + t_2 t_3))}{(-1 + t_1 t_2 t_3) (1 + t_1 t_2 t_3)}, \\
b &= \frac{2 t_2 (1 + t_2 t_3 (1 + t_3 t_1))}{(-1 + t_1 t_2 t_3) (1 + t_1 t_2 t_3)}, \\
c &= \frac{2 t_3 (1 + t_3 t_1 (1 + t_1 t_2))}{(-1 + t_1 t_2 t_3) (1 + t_1 t_2 t_3)}.
\end{align*}
We have noted that the rank jumps if $t_3(t_3-t_2)$ is a perfect square
(and, cyclicly, if $t_1(t_1-t_3)$ is a perfect square or if $t_2(t_2-t_1)$ is a perfect square).
Indeed, if we insert
$$ x= -\frac{4(t_2^2t_3-t_3+t_2) (t_3t_1^2t_2+1+t_3t_1) (t_2t_3+t_2t_3^2t_1+1)}
{t_3(-1+t_1t_2t_3)^2 (1+t_1t_2t_3)^2} $$
(note that $x+ab=\frac{b(c-b)}{t_2t_3}$)
into the equation (\ref{e2}),
we obtain
{\small
$$ y^2 =
-\frac{64(1+t_3 t_1)^2 (t_1t_2t_3-t_2-t_2^2t_3+t_3)^2 (t_2t_3+t_2t_3^2t_1+1)^2 (1+t_2 t_3)^2
(t_3 t_1^2 t_2+1+t_3 t_1)^2 (t_2-t_3)}{t_3^3(-1+t_1 t_2 t_3)^6 (1+t_1 t_2 t_3)^6}, $$ }%
which leads to the condition that $t_3(t_3-t_2)$ is a perfect square.

Thus, if we find a triple $(t_1,t_2,t_3)$ for rationals such that
\begin{equation} \label{t1t2t3}
 t_3(t_3-t_2), \quad t_1(t_1-t_3), \quad t_2(t_2-t_1)
\end{equation}
are all perfect squares, we may expect that our curve will have rank $\geq 5$
(since we started with rank $\geq 2$).

One way to satisfy conditions (\ref{t1t2t3}) is through so called almost perfect cuboids.
Indeed, if we put
$$ t_3=s_3^2,\,\,\, t_1=-s_1^2,\,\,\, t_2=s_2^2,\,\,\, s_3^2 - s_2^2 = s_4^2, $$
then we have
\begin{equation} \label{s1s2s4}
 s_1^2+s_2^2=\Box,\,\,\, s_2^2+s_4^2=\Box,\,\,\, s_1^2+s_2^2+s_4^2=\Box.
\end{equation}
Thus we get an almost perfect cuboid (only one diagonal is not an integer).
In \cite{Luijk}, one can find a parametric solution of (\ref{s1s2s4}):
\begin{align*}
s_1 &= 2(m^2+m+1)(m^2-1)^2(m^2+1+4m),\\
s_2 &= 4(m^2+m+1)(2m+1)(m^2-1)(2m+m^2),\\
s_4 &= (2m+1)(2m+m^2)(3m^2+2m+1)(m^2+2m+3).
\end{align*}
which gives
\begin{align*}
t_1 &=-4(m^2+m+1)^2 (m^2-1)^4 (m^2+1+4 m)^2,\\
t_2 &=16(m^2+m+1)^2 (2m+1)^2 (m^2-1)^2 (2m+m^2)^2,\\
t_3 &=m^2(2m+1)^2 (m+2)^2 (5m^2+8m+5)^2 (m^2+1)^2.
\end{align*}

We now present another approach which yields a two-parametric solution,
more appropriate for numerical experiments for finding specializations with higher rank.
We satisfy the first two conditions by putting
$$ t_3(t_3-t_2)=(t_3+u)^2, \quad t_1(t_1-t_3)=(t_1+v)^2 $$
and we get
$$ t_2 = -\frac{u(2t_3+u)}{t_3}, \quad t_3 = -\frac{v(2t_1+v)}{t_1}. $$
By inserting this into the third condition
$t_2(t_2-t_1)=\Box$, we get
\begin{equation} \label{uvt1}
(8uv^2-2u^2v)t_1^3+(-8u^3v+15u^2v^2+u^4+8uv^3)t_1^2+(-4u^3v^2+2v^4u+16u^2v^3)t_1+4v^4u^2=\Box.
\end{equation}
The equation (\ref{uvt1})
can be viewed as an elliptic curve over $\mathbb{Q}(u,v)$,
with an obvious point $P=[0,2u^2v^2]$.
By taking the point $2P$, we obtain
$$ t_1 = \frac{v^2 (-v+16u)}{8u (-4 v+u)}, $$
which gives
\begin{align*}
a &= -\frac{v^2 (-v+16 u) (16v^2-64u^2-v^4+16uv^3-4v^5u+v^4u^2)}{u(2+v)(4-2v+v^2)(v-2)(v^2+2v+4)(2u-v)(2u+v)(-4v+u)}, \\
b &= \frac{16u(-4v+u)v(4v-64u+16uv^2-4u^2v-v^5+4u^2v^3)}{(2+v)(4-2v+v^2)(v-2)(v^2+2v+4)(2u-v)(2u+v)(-v+16u)}, \\
c &=\frac{4(256uv-64u^2-16v^4+64u^2v^2+v^6-16v^5u)(2u-v)(2u+v)}{u(2+v)(4-2v+v^2)(v-2)(v^2+2v+4)(-v+16u)(-4v+u)}.
\end{align*}
This gives the elliptic curve with rank $\geq 5$ over $\mathbb{Q}(u,v)$.
Indeed, if we write the curve in the form $y^2 = x^3 + A x^2 + B x$, where
{\small
\begin{align*}
A &= v(256v^{13}-32v^{15}+v^{17}+140288v^9u^2+741888v^7u^4-4096v^{10}u-1167360v^8u^3 \\
&\hspace*{0.3cm}\mbox{}-21258240v^6u^5-7936v^{12}u+664832v^{10}u^3+11440128v^8u^5+32192v^{11}u^2 \\
&\hspace*{0.3cm}\mbox{}-2785824v^9u^4-32380416v^7u^6+28747776v^5u^6+6463488v^6u^7+71860224u^7v^4 \\
&\hspace*{0.3cm}\mbox{}-2205696u^8v^5+1536v^{14}u-24192v^{13}u^2-22528v^{12}u^3+591360v^{11}u^4 \\
&\hspace*{0.3cm}\mbox{}-3244800u^5v^{10}-128483328v^3u^8-12979200v^8u^7+7816v^{15}u^2-36160v^{14}u^3 \\
&\hspace*{0.3cm}\mbox{}-8616v^{13}u^4+100992v^{12}u^5-128v^{16}u-2023776v^{11}u^6+4v^{18}u-449v^{17}u^2 \\
&\hspace*{0.3cm}\mbox{}+7824v^{16}u^3-31368v^{15}u^4+2860032v^{10}u^7+70176v^{14}u^5+112296v^13u^6 \\
&\hspace*{0.3cm}\mbox{}+9461760v^7u^8-2785824v^9u^8-332160v^{12}u^7+128188416v^2u^9-37027840v^4u^9 \\
&\hspace*{0.3cm}\mbox{}-1441792v^6u^9+2659328v^8u^9+46368v^{11}u^8-6193152u^{10}v^5+515072u^{10}v^7 \\
&\hspace*{0.3cm}\mbox{}-291840u^9v^{10}+16818240v^9u^6-29425664vu^{10}+32014336u^{10}v^3+140288v^9u^{10} \\
&\hspace*{0.3cm}\mbox{}-2097152u^{11}v^2+1572864u^{11}v^4-507904u^{11}v^6-16384u^{11}v^8+65536u^{12}v^5 \\
&\hspace*{0.3cm}\mbox{}+65536u^{12}v-131072u^{12}v^3+1048576u^{11}),
\\
B &= 4(8vu^2-8u^2+16vu-v^2u+v^2+2v^3) (8vu^2+8u^2-16vu-v^2u-v^2+2v^3) \\
 &\hspace*{0.3cm}\mbox{}\times (-16v^2+64u^2+v^4-16v^3u) (4v-64u+16v^2u-4vu^2-v^5+4v^3u^2) \\
 &\hspace*{0.3cm}\mbox{}\times (2vu^2-16u^2+2vu+8v^2u-4v^2-v^3) (2vu^2+16u^2-2vu+8v^2u+4v^2-v^3) \\
 &\hspace*{0.3cm}\mbox{}\times (16vu-4u^2-v^4+4v^2u^2) (16v^2-64u^2-v^4+16v^3u-4v^5u+v^4u^2) \\
 &\hspace*{0.3cm}\mbox{}\times (-v+16u)^2 (-4 v+u)^2 u^2 v^3,
\end{align*}}%
then five independent points of infinite order are
{\small
\begin{align*}
P &=[-4(4v-64u+16v^2u-4vu^2-v^5+4v^3u^2) (16v^2-64u^2-v^4+16v^3u-4v^5u+v^4u^2) \\
&\hspace*{0.3cm}\mbox{}\times (-4v+u)^2 u^2 (-v+16 u)^2 v^3, \\
&8(64v^2u^2-64u^2-16v^5u+256vu+v^6-16v^4) (4v-64u+16v^2u-4vu^2-v^5+4v^3u^2) \\
&\hspace*{0.3cm}\mbox{}\times(16v^2-64u^2-v^4+16v^3u-4v^5u+v^4u^2) (2u-v)^2 (2u+v)^2 (-4v+u)^2 u^2 (-v+16u)^2 v^3],
\\
R &=[4(16vu-4u^2-v^4+4v^2u^2) (8vu^2+8u^2-16vu-v^2u-v^2+2v^3) v (-v+16 u) (-4v+u) u \\
&\hspace*{0.3cm}\mbox{}\times(8vu^2-8u^2+16vu-v^2u+v^2+2v^3) (16v^2-64u^2-v^4+16v^3u-4v^5u+v^4u^2), \\
&4(8vu^2-8u^2+16vu-v^2u+v^2+2v^3) (16v^2-64u^2-v^4+16v^3u-4v^5u+v^4u^2) \\
&\hspace*{0.3cm}\mbox{}\times (8vu^2+8u^2-16vu-v^2u-v^2+2v^3) (16vu-4u^2-v^4+4v^2u^2)  \\ &\hspace*{0.3cm}\mbox{}\times (8vu^2+8u^2+32vu-16v^2u-4v^2-v^3) (8vu^2-8u^2-32vu-16v^2u+4v^2-v^3) \\
&\hspace*{0.3cm}\mbox{}\times (2u+v) (2u-v) v^2 (-v+16u) (-4v+u) u],
\\
T_1 &=[16(16vu-4u^2-v^4+4v^2u^2) (2vu^2-16u^2+2vu+8v^2u-4v^2-v^3) (-v+16u) (-4v+u) u \\
&\hspace*{0.3cm}\mbox{}\times(2vu^2+16u^2-2vu+8v^2u+4v^2-v^3) (4v-64u+16v^2u-4vu^2-v^5+4v^3u^2), \\
&8(16vu-4u^2-v^4+4v^2u^2) (2vu^2-16u^2+2vu+8v^2u-4v^2-v^3) (-v+16u) (-4v+u) u \\
&\hspace*{0.3cm}\mbox{}\times(2vu^2+16u^2-2vu+8v^2u+4v^2-v^3) (-v+16u-4v^2u+vu^2) (8u^2-vu+2v^2) \\ &\hspace*{0.3cm}\mbox{}\times(v^6-16v^5u+256vu-64u^2-16v^4+64v^2u^2) (4v-64u+16v^2u-4vu^2-v^5+4v^3u^2)],
\\
T_2 &=[-4(8vu^2-8u^2+16vu-v^2u+v^2+2v^3) (-16v^2+64u^2+v^4-16v^3u) \\
&\hspace*{0.3cm}\mbox{}\times (16vu-4u^2-v^4+4v^2u^2) (8vu^2+8u^2-16vu-v^2u-v^2+2v^3) \\ &\hspace*{0.3cm}\mbox{}\times (16v^2-64u^2-v^4+16v^3u-4v^5u+v^4u^2) (-4v+u) u/v^2, \\
&4(8vu^2-8u^2+16vu-v^2u+v^2+2v^3) (-16v^2+64u^2+v^4-16v^3u) \\
&\hspace*{0.3cm}\mbox{}\times (16vu-4u^2-v^4+4v^2u^2) (8vu^2+8u^2-16vu-v^2u-v^2+2v^3) (2u+v) (2u-v) \\ &\hspace*{0.3cm}\mbox{}\times (8u^2-16vu-v^2) (-16v^4+64v^2u^2+v^6-16v^5u+256vu-64u^2) \\ &\hspace*{0.3cm}\mbox{}\times (16v^2-64u^2-v^4+16v^3u-4v^5u+v^4u^2) (-4v+u) u/v^3],
\\
T_3 &=[(-v+16u) (4v-64u+16v^2u-4vu^2-v^5+4v^3u^2) (16v^2-64u^2-v^4+16v^3u-4v^5u+v^4u^2) \\ &\hspace*{0.3cm}\mbox{}\times (-16v^2+64u^2+v^4-16v^3u)(2u^2+8vu-v^2)^2, \\
&2(-16v^2+64u^2+v^4-16v^3u) (2u-v) (2u+v) (8vu^2-8u^2-32vu-16v^2u+4v^2-v^3) \\
&\hspace*{0.3cm}\mbox{}\times (8vu^2+8u^2+32vu-16v^2u-4v^2-v^3) (-v+16u-4v^2u+vu^2) (2u^2+8vu-v^2) \\ &\hspace*{0.3cm}\mbox{}\times(16v^2-64u^2-v^4+16v^3u-4v^5u+v^4u^2) (4v-64u+16v^2u-4vu^2-v^5+4v^3u^2) (-v+16u)].
\end{align*}}%
Here the point $P$ corresponds to $[0,abc]$ on $y^2=(x+ab)(x+ac)(x+bc)$,
the point $R$ satisfies $2R=S$, where $S$ corresponds to $[1,rst]$ on $y^2=(x+ab)(x+ac)(x+bc)$,
the point $T_1$ corresponds to the condition $t_3(t_3-t_2) = \Box$,
the point $T_2$ corresponds to the condition $t_1(t_1-t_3) = \Box$, while
the point $T_3$ corresponds to the condition $t_2(t_2-t_1) = \Box$.
Since the specialization map in a homomorphism, it suffices to find a specialization
$(u_0,v_0)$ for which the points $P$, $R$, $T_1$, $T_2$ and $T_3$ are independent points of
infinite order on $y^2 = x^3 + A x^2 + B x$. We checked that this is the case for $(u_0,v_0)=(2,1)$,
since the points
$[170605, 39532697]$, $[302665, -66247363]$, $[795565, -637321303]$,  $[-447095, -24260803]$,
$[8673115/4, -25165674989/8]$ are independent on $y^2 = x^3 + 21361758597 x^2 -28803989016278714304 x$.


\medskip

Now we search for specializations $(u,v)$ with higher rank,
in particular with rank $11$ and $12$.
We use a sieving methods similar to those used e.g. in \cite{ADP,D-Peral-RACSAM}.
We searched for curves with
relatively large Mestre-Nagao sum
$$ S(N,E)=\sum_{p=2}^N\dfrac{-a_p+2}{p+1-a_p}\,\log p, $$
where $a_p=a_p(E)=p+1- \# E(\mathbb{F}_p)$, since it is experimentally known \cite{M1,N} that
we may expect that high rank curves have large $S(N,E)$, and large Selmer rank
(as implemented in \texttt{mwrank} with option \texttt{-s}).
In search for rank $12$ curves we also use the condition that the root-number is equal to $1$
(conjecturally this implies that rank is even).
We searched also in some restricted subfamilies, including e.g. $u=v$.
We implemented the sieving algorithm in {\tt Pari} \cite{Pari}.
For the curves which pass our searching conditions, we calculate the rank
by Cremona's program \texttt{mwrank} \cite{Mwrank}.

We find curves with rank $11$ for the following parameters:
$(u,v)=$
\begin{align*}
& \left(\frac{11}{24}, \frac{5}{9}\right), \left(-\frac{145}{6},  \frac{29}{12}\right),
\left(\frac{136}{19}, \frac{68}{5}\right),
\left(-\frac{16}{77}, \frac{4}{21}\right), \left(\frac{473}{705}, \frac{43}{47}\right),
\left(-\frac{89}{135}, \frac{89}{45}\right), \\
& \left(-\frac{62}{43}, \frac{93}{43}\right),
\left(\frac{71}{273}, \frac{142}{91}\right), \left(\frac{224}{67}, \frac{7}{2}\right),
\left(-\frac{1032}{923}, \frac{172}{71}\right),
\left(-\frac{1501}{87}, \frac{158}{87}\right), \\
& \left(\frac{1358}{1007}, \frac{194}{53}\right),
\left(-\frac{2072}{1819}, \frac{148}{107}\right),
\left(\frac{454}{481}, \frac{227}{37}\right), \left(\frac{77}{173}, \frac{77}{173}\right),
\left(\frac{163}{137}, \frac{163}{137}\right).
\end{align*}
The details (minimal Weierstrass equation, torsion points and independent points of infinite order)
are given in \cite{D-tors}.
Let us mention that the curve corresponding to $(u,v)=(-\frac{62}{43}, \frac{93}{43})$,
i.e.
$$ \{ a,b,c \} = \left\{ \frac{21409906185}{74591676404}, -\frac{31580198976}{18647919101}, -\frac{10309975195}{18647919101} \right\}, $$
with the minimal Weierstrass equation
\begin{align*}
y^2 + xy &= x^3 - x^2 - 21252276640652798739707819217x \\
&\hspace*{0.3cm} \mbox{}+ 938627524108684110053910801619511357084941,
\end{align*}
has the minimal discriminant among all known curves with rank $11$ and torsion group
$\mathbb{Z}/2\mathbb{Z} \times \mathbb{Z}/2\mathbb{Z}$.

Finally, we found a curve with rank $12$ for $(u,v)=(-\frac{95}{33}, \frac{50}{57})$, i.e.
$$ \{a,b,c\}=\left\{ \frac{6125241375}{11907531272}, \frac{5535371271425}{14277129995128},  -\frac{273138178560}{153430695649} \right\}, $$
with the minimal Weierstrass equation
\begin{align*}
y^2 + xy + y &= x^3 - x^2 - 1444491707528591356856089186460491195711268950880x \\
 &\hspace*{0.3cm} \mbox{}+ 559921583779625421248683584939561762456224290170437461555851482041439747,
\end{align*}
the torsion points
\begin{align*}
&\mathcal{O},\,\, [910954389920845836020349, -455477194960422918010175], \\
&[-5448727291190824028230629/4, 5448727291190824028230625/8], \\
&[451227432876860171037309, -225613716438430085518655],
\end{align*}
and $12$ independent points of infinite order
\begin{align*}
P_1 &= [158850932500649609134809, 578334775816714524616276221704042845], \\
P_2 &= [351104017200784386392209, 309897966944945116194624198332593845], \\
P_3 &= [-427722660290928813983135, -1048576645526111528109185629948786727], \\
P_4 &= [954500781939375762742909, 225326008863345220543071618783370945], \\
P_5 &= [423679598259676591990909, 154829810959547852593332987635966145], \\
P_6 &= [1535808449095818094207905, 1401421444080498380369785533616999513], \\
P_7 &= [444801887422056021535383, 73569216148613399817347986859758945], \\
P_8 &= [-1206006015871044278678751, -740210245609217615143269452335454375], \\
P_9 &= [-192562292438693523617091, -911556889640548767064630159456313855], \\
P_{10} &= [10508879668527356682921249, 33851800053181168926568362825476385625], \\
P_{11} &= [951514410733369555670349, 216676520921276805299703311439049825], \\
P_{12} &= [-7355680099955426717481581/81, -605705671933225602690651446390633849125/729].
\end{align*}

Let us also mention a minor, somewhat related result:
for $t_1=\frac{44}{29}$, $t_2=\frac{17}{42}$, $t_3=\frac{3}{44}$, i.e.
$a=\frac{815848}{164547}$, $b=\frac{1512524}{1810017}$, $c=\frac{32060}{201113}$,
we get the elliptic curve
$$ y^2 = x^3 + x^2 - 193936360896469946772176x + 29453641253718130506136229522416740 $$
with rank $10$,
which is the curve with smallest known conductor
among curves with rank $10$ and torsion group $\mathbb{Z}/2\mathbb{Z} \times \mathbb{Z}/2\mathbb{Z}$.
It is obtained by brute-force search (not in parametric families)
within triples with small $t_1$, $t_2$, $t_3$.

\section{An infinite family of elliptic curves with rank $\geq 7$}

The construction of the two-parametric family of curves with rank $\geq 5$
from the previous section is related with the construction
from our joint paper with Juli\'an Aguirre \cite{ADP}.
In \cite{ADP}, we constructed a two-parametric family of curves with rank $\geq 4$
over $\mathbb{Q}(m,n)$, and by choosing $n=7/3$ we obtained a family with rank $\geq 5$ over $\mathbb{Q}(m)$.
It can be checked that by taking $m = -\frac{20(4u^2-1)}{9u(u+4)}$,
we obtain the same family as the family
obtained from our new two-parametric family by specializing $v=-1$.

It is shown in \cite{D-Peral-RACSAM} that inserting $n=7/3$ already in the family from \cite{ADP}
with rank $\geq 3$ over $\mathbb{Q}(a,n)$, gives a simple family with rank $\geq 4$ over $\mathbb{Q}(a)$,
which is very suitable for constructing subfamilies with higher rank.
That family is
$$ y^2 = x^3 + A(a)x^2 + B(a)x, $$
where
   \begin{align*}
A(a)&=-2 (-51200 + 109440 a + 38880 a^2 + 55404 a^3 + 6561 a^4),\\
B(a)&=243 a^2 (20 + 3 a) (-4 + 9 a) (16 + 9 a) (80 + 9 a) (320 + 81 a^2),
\end{align*}
and the  $x$-coordinates of four independent points of infinite order are
\begin{equation}\label{xcoord}
\begin{aligned}
x_1&= 81 a^2 (-4 + 9 a) (80 + 9 a),\\
x_2&=27 a (20 + 3 a) (-4 + 9 a) (80 + 9 a),\\
x_3&=\frac{1}{441} (-4 + 9 a) (80 + 9 a) (160 + 171 a)^2,\\
x_4&=3 (20 + 3 a) (-4 + 9 a) (320 + 81 a^2).
\end{aligned}
\end{equation}
There are several substitutions which give subfamilies with rank $\geq 6$:
\begin{align*}
 a &= -\frac{2 (-27 + 13 w_1^2) (-13 + 27 w_1^2)}{9 (9 + 178 w_1^2 + 9 w_1^4)}, \\ 
 a &= -\frac{ 64 (831744 - 40128 w_2 + 4288 w_2^2 - 44 w_2^3 + w_2^4)}
 {9 (-1520 + 88 w_2 + w_2^2) (-2736 - 264 w_2 + 5 w_2^2)}, \\ 
 a &= \frac{10732176 + 628992 w_3 + 19192 w_3^2 + 576 w_3^3 + 9 w_3^4}{36 w_3 (27 + w_3) (364 + 9 w_3)}, \\ 
 a &= \frac{5 (-10 + 6 w_4 + w_4^2) (-18 - 18 w_4 + 5 w_4^2)}{9 (12 - 2 w_4 + w_4^2) (3 - w_4 + w_4^2)}, \\ 
 a &= \frac{ 5 (584820 + 135432 w_5 - 18288 w_5^2 + 396 w_5^3 + 5 w_5^4)}{9 (684 - 66 w_5 + w_5^2) (171 - 33 w_5 + w_5^2)}. \\ 
\end{align*}
The first four substitutions were already given in \cite{D-Peral-RACSAM}, while the fifth
substitution is new.

In order to find infinite families with rank $\geq 7$,
we try to find intersections of these five families with rank $\geq 6$.
We compare their $j$-invariants by factorizing their difference
and seeking for the factors which correspond to curves with genus $1$.

If we compare the second and third substitution, we find two suitable factors,
which give the following conditions:
\begin{align}
& w_2^2 w_3^2+72w_2^2w_3+88w_2w_3^2+1820w_2^2-1520w_3^2-96096w_2-65664w_3-995904 =0, \label{t2t3a}\\
& 5w_2^2w_3^2+21w_2^2w_3-264w_2w_3^2+3276w_2^2-2736w_3^2+288288w_2-196992w_3-4979520 =0. \label{t2t3b}
\end{align}
Both conditions lead to
\begin{equation} \label{quarticw2w3}
 54w_3^4+2736w_3^3+66592w_3^2+2987712w_3+64393056 = \Box.
\end{equation}
This quartic is birationally equivalent to the elliptic curve
$$ y^2 = x^3+x^2-28174550x+45644288448 $$
with rank equal to $3$, hence the elliptic curve, and also the quartic, have
infinitely many rational solutions. Many of them produce curves with rank $= 7$,
e.g.  $t_3=-234, -30, -18, 26, 42, 94, -\frac{202}{3}, -\frac{182}{3}, -\frac{14}{3}$.

Consider the four points given by (\ref{xcoord}) and additional two points
corresponding to the second and third substitutions. The second substitution gives the curve
\begin{equation} \label{62}
y^2 = x^3 + a_{62} x^2 + b_{62} x,
\end{equation}
where
{\small
\begin{align*}
a_{62} &= 79573w_2^{16}+2281840w_2^{15}-791687936w_2^{14}-34844285696w_2^{13}+3065917324288w_2^{12} \\
&\hspace*{0.3cm}\mbox{}+556971294060544w_2^{11}-64165839736733696w_2^{10}+3360211454234263552w_2^9-130403990149389221888w_2^8 \\
&\hspace*{0.3cm}\mbox{}+3064512846261648359424w_2^7-53369552205989831245824w_2^6+422490869190468915167232w_2^5\\ &\hspace*{0.3cm}\mbox{}+2120995723090424777146368w_2^4-21983951517250398896259072w_2^3-455536370311599498486349824w_2^2 \\
&\hspace*{0.3cm}\mbox{}+1197427029434259336824094720w_2+38082411231292796255084740608, \\
b_{62} &= -5184(w_2^4-44w_2^3+4288w_2^2-40128w_2+831744)^2 (w_2^4+352w_2^3-50720w_2^2+321024w_2+831744) \\ &\hspace*{0.3cm}\mbox{}\times (3w_2^4+352w_2^3+15328w_2^2-642048w_2+5822208) (7w_2^4-704w_2^3+15328w_2^2+321024w_2+2495232) \\
&\hspace*{0.3cm}\mbox{}\times (7w_2^4-176w_2^3+11680w_2^2-160512w_2+5822208) (7w_2^4+352w_2^3-61664w_2^2+321024w_2+5822208) \\
&\hspace*{0.3cm}\mbox{}\times (59w_2^4+3344w_2^3-572128w_2^2+3049728w_2+49072896),
\end{align*}}%
and six independent points of infinite order with $x$-coordinates:
{\small
\begin{align*}
x_{21} &= -576(w_2^4-44w_2^3+4288w_2^2-40128w_2+831744)^2 (7w_2^4-176w_2^3+11680w_2^2-160512w_2+5822208) \\
&\hspace*{0.3cm}\mbox{}\times (7w_2^4+352w_2^3-61664w_2^2+321024w2+5822208), \\
x_{22} &= 36(w_2^4-44w_2^3+4288w_2^2-40128w_2+831744) (7w_2^4-176w_2^3+11680w_2^2-160512w_2+5822208) \\ &\hspace*{0.3cm}\mbox{}\times (7w_2^4+352w_2^3-61664w_2^2+321024w_2+5822208) (59w_2^4+3344w_2^3-572128w_2^2+3049728w_2+49072896), \\
x_{23} &= -16/49(7w_2^4-176w_2^3+11680w_2^2-160512w_2+5822208) (7w_2^4+352w_2^3-61664w_2^2+321024w_2+5822208) \\ &\hspace*{0.3cm}\mbox{}\times (13w_2^4-2552w_2^3+330784w_2^2-2327424w_2+10812672)^2, \\
x_{24} &= -27/4(3w_2^4+352w_2^3+15328w_2^2-642048w_2+5822208) (7w_2^4-704w_2^3+15328w_2^2+321024w_2+2495232) \\ &\hspace*{0.3cm}\mbox{}\times (7w_2^4-176w_2^3+11680w_2^2-160512w_2+5822208) (59w_2^4+3344w_2^3-572128w_2^2+3049728w_2+49072896), \\
x_{25} &= -108(w_2^2-912)^2 (w_2^4+352w_2^3-50720w_2^2+321024w_2+831744) \\ &\hspace*{0.3cm}\mbox{}\times (7w_2^4+352w_2^3-61664w_2^2+321024w_2+5822208) (59_w2^4+3344w_2^3-572128_w2^2+3049728w_2+49072896), \\
x_{26} &= 324(w_2^2-912)^2 (w_2^4+352w_2^3-50720w_2^2+321024w_2+831744) \\ &\hspace*{0.3cm}\mbox{}\times (7w_2^4-176w_2^3+11680w_2^2-160512w_2+5822208) (7w_2^4+352w_2^3-61664w_2^2+321024w_2+5822208).
\end{align*}}%
The third condition gives the curve
\begin{equation} \label{63}
 y^2 = x^3 + a_{63} x^2 + b_{63} x,
\end{equation}
where
{\small
\begin{align*}
a_{63} &= -13122w_3^{16}-7348320w_3^{15}-1570137696w_3^{14}-206172584064w_3^{13}-19541430237312w_3^{12} \\
&\hspace*{0.3cm}\mbox{}-1402008391816704w_3^{11}-77606011598363136w_3^{10}-3410103604914358272w_3^9-123219415654113963008w_3^8 \\
&\hspace*{0.3cm}\mbox{}-3723833136566479233024w_3^7-92542375014630498607104w_3^6-1825654232153731017572352w_3^5 \\
&\hspace*{0.3cm}\mbox{}-27787335201034030779236352w_3^4-320143070559304939026382848w_3^3-2662401630093588063697895424w_3^2 \\
&\hspace*{0.3cm}\mbox{}-13606503227295711027839631360w_3-26532681293226636504287281152, \\
b_{63} &= 81(w_3^4+72w_3^3+8504w_3^2+550368w_3+10732176) (3w_3^4+144_w3^3+3160_w3^2+157248_w3+3577392) \\ &\hspace*{0.3cm}\mbox{}\times (3w_3^4+1152w_3^3+71144w_3^2+1257984w_3+3577392) (9w_3^4+504w_3^3+8504w_3^2+78624w_3+1192464) \\
&\hspace*{0.3cm}\mbox{}\times (9w_3^4+576w_3^3+19192w_3^2+628992w_3+10732176)^2 (9w_3^4+1152w_3^3+58040w_3^2+1257984w_3+10732176) \\
&\hspace*{0.3cm}\mbox{}\times (9w_3^4+2736w_3^3+164872w_3^2+2987712w_3+10732176),
\end{align*}}%
and six independent points of infinite order with $x$-coordinates:
{\small
\begin{align*}
x_{31} &= 9(3w_3^4+144w_3^3+3160w_3^2+157248w_3+3577392) (3w_3^4+1152w_3^3+71144w_3^2+1257984w_3+3577392) \\ &\hspace*{0.3cm}\mbox{}\times (9w_3^4+576w_3^3+19192w_3^2+628992w_3+10732176)^2, \\
x_{32} &= 9(3w_3^4+144w_3^3+3160w_3^2+157248w_3+3577392) (3w_3^4+1152w_3^3+71144w_3^2+1257984w_3+3577392) \\ &\hspace*{0.3cm}\mbox{}\times (9w_3^4+576w_3^3+19192w_3^2+628992w_3+10732176) (9w_3^4+2736w_3^3+164872w_3^2+2987712w_3+10732176), \\
x_{33} &= 1/49(3w_3^4+144w_3^3+3160w_3^2+157248w_3+3577392) (3w_3^4+1152w_3^3+71144w_3^2+1257984w_3+3577392) \\ &\hspace*{0.3cm}\mbox{}\times (171w_3^4+16704w_3^3+753128w_3^2+18240768w_3+203911344)^2, \\
x_{34} &= 27(w_3^4+72w_3^3+8504w_3^2+550368w_3+10732176) (3w_3^4+144w_3^3+3160w_3^2+157248w_3+3577392) \\ &\hspace*{0.3cm}\mbox{}\times (9w_3^4+504w_3^3+8504w_3^2+78624_w3+1192464) (9w_3^4+2736w_3^3+164872w_3^2+2987712w_3+10732176), \\
x_{35} &= 27(w_3^2-1092)^2 (3w_3^4+1152w_3^3+71144w_3^2+1257984w_3+3577392) \\
&\hspace*{0.3cm}\mbox{}\times (9w_3^4+1152w_3^3+58040w_3^2+1257984w_3+10732176) (9w_3^4+2736w_3^3+164872w_3^2+2987712w_3+10732176), \\
x_{36} &= 81(w_3^2+54w_3+1092)^2 (3w_3^4+144w_3^3+3160w_3^2+157248w_3+3577392) \\
&\hspace*{0.3cm}\mbox{}\times (3w_3^4+1152w_3^3+71144w_3^2+1257984w_3+3577392) (9w_3^4+1152w_3^3+58040w_3^2+1257984w_3+10732176).
\end{align*}}%

By factorizing the expression $a_{62}^2 b_{63} - a_{63}^2 b_{62}$ we see that
for pairs $(w_2,w_3)$ satisfying the conditions (\ref{t2t3a}) or (\ref{t2t3b}) it holds
$b_{63}/a_{63}^2 = b_{62}/a_{62}^2$. Hence, the curves (\ref{62}) and (\ref{63}) are
birationally equivalent.
In the same way we check that for such pair $(w_2,w_3)$ it holds
$x_{31}/a_{63} = x_{21}/a_{62}$,  $x_{32}/a_{63} = x_{22}/a_{62}$, $x_{33}/a_{63} = x_{23}/a_{62}$,
$x_{34}/a_{63} = x_{24}/a_{62}$, $x_{35}/a_{63} = x_{25}/a_{62}$, while
$x_{36}/a_{63} \neq x_{26}/a_{62}$.
Therefore, we have seven points on (\ref{63}) with $x$-coordinates
\begin{equation} \label{7points}
x_{31},\,\, x_{32},\,\, x_{33},\,\, x_{34},\,\, x_{35},\,\, x_{36},\,\,x_{26} a_{63}/a_{62}.
\end{equation}
By taking the specialization $(w_2,w_3)=(-\frac{76}{3},26)$ we obtain the curve
\begin{align*}
y^2 &=x^3-163531808801344950045916528640000x^2 \\ &\hspace*{0.3cm}\mbox{}+6680706316011654681276493655189069731350803361465165152256000000x
\end{align*}
and we checked that seven corresponding points with $x$-coordinates
\begin{align*}
& 38540677847903454008558223360000,\,\, 178922409809838644555667210240000, \\
& 72051389475320867247399895040000,\,\, 66579605091474988619076835737600, \\
& 13362426543070313045072805888000,\,\, 126710845595682509808491456102400 \\
& \hspace*{3.0cm}\mbox{and} \,\,\,\, 2179385680764224839490312601600
\end{align*}
are independent points of infinite order on this curve.
By Silverman's specialization theorem \cite[Theorem III.11.4]{Silverman2}, we conclude that seven points (\ref{7points})
are independent points on (\ref{63}) for infinitely many rational values of $w_3$
satisfying the quartic equation (\ref{quarticw2w3}) and the corresponding values $w_2$ satisfying
(\ref{t2t3a}) or (\ref{t2t3b}).
Thus we proved that there are indeed infinitely many elliptic curves induced by rational
Diophantine triples with rank $\geq 7$.

\medskip

Analogous result can be obtained by considering the second and fifth substitution for the parameter $a$.
Here the conditions are
\begin{align*}
&9w_2^2w_5-4w_2w_5^2-198w_2^2+528w_5^2+1368w_2-8208w_5 = 0, \\
&11w_2^2w_5^2-171w_2^2w_5-76w_2w_5^2+25992w_2+155952w_5-3430944 = 0,
\end{align*}
and they lead to the quartic
$$
w_5^4-1188w_5^3+43920w_5^2-406296w_5+116964 = \Box,
$$
which is equivalent to the elliptic curve
$$ y^2 = x^3-x^2-124056x-10126800 $$
with rank equal to $1$, so we again have infinitely many rational solutions.
These solutions give seven points on the curve (\ref{62}).
By taking the specialization $(w_2,w_5)=(\frac{6392}{99},\frac{6392}{99})$
we can check that these seven points are indeed independent,
and by Silverman's specialization theorem we conclude that we obtained another infinite family
of curves with rank $\geq 7$.

\bigskip

{\bf Acknowledgements.}
The authors would like to thank Luka Lasi\'c for useful comments on the previous version of this paper.
A.D. was supported by the Croatian Science Foundation under the project no.~IP-2018-01-1313.
He also acknowledges support from the QuantiXLie Center of Excellence, a project
co-financed by the Croatian Government and European Union through the
European Regional Development Fund - the Competitiveness and Cohesion
Operational Programme (Grant KK.01.1.1.01.0004).
J. C. P. was supported by the grant: MTM2014-52347-C2-1-R.

\end{document}